\documentclass[11pt]{article}
\usepackage{a4wide}
\usepackage{amsmath}

\usepackage{amsfonts}
\usepackage{amssymb}
\usepackage{euscript}
\newenvironment{proof}{\noindent {\it Proof.~~}\ }{\  \rule{1mm}{2mm}\medskip}

\newenvironment{proof*}{\noindent {\it Proof.~~}\ }{}

\newenvironment{proofof}[2]{\noindent {\it Proof of #1}~#2: \
}{~\rule{1mm}{2mm}\medskip}
\newtheorem{theorem}{Theorem}
\newtheorem{lemma}[theorem]{Lemma}
\newtheorem{corollary}[theorem]{Corollary}

\newtheorem{conj}[theorem]{Conjecture}

\def\Z{\mathbb Z}

\begin{document}
\author{ Yahya O. Hamidoune\thanks{Universit\'e Pierre et Marie Curie,    Paris {\tt yha@ccr.jussieu.fr} }
}
\title{A weighted generalization of  Gao's n+D-1 Theorem}
\maketitle

\begin{abstract}

Let $G$ denotes  a finite abelian group of order $n$ and Davenport constant $D$,
and put $m= n+D-1$.
 Let $x=(x_1, \cdots ,x_m)\in G^m$ be a sequence with a maximal repetition $\ell$
 attained by $x_m$ and put $r=\min(D,\ell)$.
 Let  $w=(w_1, \cdots ,w_{m-r})\in \Z^{m-r}.$

Then there are an  $n$-subset $I\subset [1,m-r]$ and an injection
$f: I\mapsto  [1,m]$, such that $m\in f(I)$ and
  $$\sum_{ i\in I}w_{i}x_{f({i})}=(\sum_{ i\in I}w_{i})x_{m}.$$
\end{abstract}

\section{Introduction }

Let $G$ be an abelian group of order $n$. The {\em Davenport constant} of $G$, denoted by $D(G)$
is the maximal integer $k$ such that there is a sequence of elements of $G$ with length $k-1$
having no  nonempty  zero-sum subsequence. The investigation  of zero-sum subsequences with length $n$ starts with a result of  Erd{\H o}s,  Ginzburg and  Ziv \cite{egz}, stating that
every sequence of elements of $G$ with length $\ge 2n-1$ has a zero-sum subsequence of length $n$.
A result of Gao \cite{gaotnd} stating that
every sequence of elements of $G$ with length $\ge n+D(G)-1$ has a  zero-sum subsequence of length $n$,
unifies these two areas of Zero-sum problems.
For more details on these questions and some of their applications, the reader may refer to the book of Geroldinger, Halter-Koch \cite{geroldinger}
and to the survey paper of Caro \cite{caro}.

Attempts were made to generalize Zero-sum problems to the weighted case. Examples of such results
 may be found in the survey paper of Caro \cite{caro}, a paper by the author \cite{string}   and two more recent papers by  Gao-Jin \cite{GJ} and  Grynkiewicz \cite{david}.

Let $x =(x_1,\cdots ,x_m)\in E^m,$ where  $E$ is a set. As usual one may consider $x$ as a map form
$[1,n]$ into $E$. We
shall put $\rho(x)=\max \{ |x^{-1}(a)| ; a\in E\}$. Notice that $\rho (x)$ is the maximal repetition of the sequence $x$.

We shall prove some  weighted generalizations of Gao's Theorem.
Our main result in this note is the following:

\begin{theorem}\label{weightedgao}

Let $G$ denotes  a finite abelian group of order $n$ and Davenport constant $D$.
 Let $x=(x_1, \cdots ,x_m)\in G^m$  and let  $w=(w_1, \cdots ,w_m)\in \Z^m.$
Assume moreover that $|x^{-1}(a)|\le \ell$, for all $a\in G$
and that  $m= n+D-\min (D,\ell)-1$.

Then there are
a subset $I\subset [1,m]$ and an injection
$f: I\mapsto  [1,m]$, such that

$n-\min (D,\ell) \le |I| \le n-1$ and
  $$\sum_{ i\in I}w_{i}x_{f(i)}=0.$$

\end{theorem}

\begin{corollary}\label{baricenter}
Let $G$ denotes  a finite abelian group of order $n$ and Davenport constant $D$,
and put $m= n+D-1$.
 Let $x=(x_1, \cdots ,x_m)\in G^m$ be such that
 $\max _{a\in G}|x^{-1}(a)|= \ell=|x^{-1}(x_m)|$ and put $r=\min(D,\ell)$.
 Let  $w=(w_1, \cdots ,w_{m-r})\in \Z^{m-r}.$

Then there are an  $n$-subset $I\subset [1,m-r]$ and an injection
$f: I\mapsto  [1,m]$, such that $m\in f(I)$ and
  $$\sum_{ i\in I}w_{i}x_{f({i})}=(\sum_{ i\in I}w_{i})x_{m}.$$
\end{corollary}

With $w_i=1$ for all $i$, Corollary \ref{baricenter} reduces to Gao $n+D-1$-Theorem.

We need the following result:

\begin{lemma}\label{word0} (folklore)
Let $G$ be a finite abelian group of order $n$.
  Let $x=(x_1, \cdots ,x_n)\in G^n$ be a sequence of  elements of $G$
  with a maximal repetition $\le k$.
  Then $x$ has a nonempty zero-sum subsequence with length $\le k$.
\end{lemma}

Lemma \ref{word0} is now a standard tool in   Zero-sum  Theory
 \cite{geroldinger}.  A proof of Lemma \ref{word0}  requires an Addition theorem. The oldest such a result implying  easily Lemma \ref{word0} is Shepherdon's Theorem
\cite{sheph}. An almost identical proof follows  by Scherck's Theorem \cite{scherck}. When $G$ is cyclic, Lemma \ref{word0} is a special case
of Conjecture 4 of  Erd\H{o}s and  Heilbronn \cite{EH}. In a note added in proofs, Erd\H{o}s and  Heilbronn \cite{EH}
mentioned that Flor proved this conjecture. A generalization of Lemma \ref{word0} to non necessarily abelian groups is proved in
\cite{hwords} using a result of Kemperman \cite{Kemp}.
\section{Weighted sums}
Let $G$ denotes  a finite abelian group of order $n$.

We need the following weighted generalization of Lemma \ref{word0}:

\begin{lemma}\label{word1}

  Let $x=(x_1, \cdots ,x_n)\in G^n$  and let  $w=(w_1, \cdots ,w_n)\in \Z^n.$
  Assume that $\max _{a\in G}|x^{-1}(a)|\le  \ell$.
  Then there are a subset $I\subset [1,n]$ with $1\le |I|\le \ell$, and an injection
  $f:I\mapsto [1,n]$ such that
  $$\sum_{ i\in I}w_{i}x_{f(i)}=0.$$
\end{lemma}
\begin{proof}
Put $r=\rho(x)$ and $s=\rho(w)$. Take $a\in G$ such that $|x^{-1}(a)|=r$, and
$b\in \Z$ such that $|w^{-1}(b)|=s$.

Assume first that $s\le r$. By Lemma \ref{word0}, applied to $w$, there is a subset $I\subset [1,n]$
with $|I|\le r$ such that $\sum_{i\in I}w_{i}\equiv 0$ mod $n$. Take an arbitrary injection $f:I\mapsto x^{-1}(a)$.
We have $\sum_{i\in I}w_ix_{f(i)}=\sum_{i\in I}w_ia=(\sum_{i\in I}w_i)a=0.$

Assume now that $s> r$. By Lemma \ref{word0}, there is a subset $J\subset [1,n]$
with $|J|\le r$ such that $\sum_{i\in J}x_{i}=0.$
Take a subset $L \subset w^{-1}(b)$ such that $|L|=|J|$
Take an arbitrary bijection $f:L\mapsto J$.
We have $\sum_{i\in L} w_ix_{f(i)}=\sum_{i\in J} bx_i=b(\sum_{i\in J}x_i)=0.$
\end{proof}

This Lemma suggests the following conjecture:

\begin{conj}\label{folkconj}.

  Let $x=(x_1, \cdots ,x_n)\in G^n$  and let  $w=(w_1, \cdots ,w_k)\in \Z^k.$
  Assume that $\max _{a\in G}|x^{-1}(a)|\le  k$.
  Then there are a nonempty subset $I\subset [1,k]$  and an injection
  $f:I\mapsto [1,n]$ such that
  $$\sum_{ i\in I}w_{i}x_{f(i)}=0.$$
\end{conj}

Conjecture \ref{folkconj} follows easily By Scherck's Theorem \cite{scherck} or by some other additive ingredients
if $\gcd (n,w_i)=1$, for all $i$. It holds trivially if $D(G)\le k.$ If true, this conjecture has very interesting
implications.

Let  $x=(x_1, \cdots ,x_m)\in G^m$  and let  $w=(w_1, \cdots ,w_m)\in \Z^m.$
Let $f$ be a function on $[1,m]$.
 We shall write $<S>_f=\sum_{ i\in S}w_ix_{f({i})}$, for a subset $S$ on which $f$ is defined.
 Let $I\subset [1,m]$ and let  $f:I\mapsto [1,m]$. We shall say that the pair $(I,f)$ is  $k$--{\em shellable} if there is a partition $I=I_1\cup \cdots \cup I_t$
such that
$1\le |I_i|\le k$ and $<I_i>_f=0$, for all $1\le i\le t$. We shall call the partition $I=I_1\cup \cdots \cup I_t$ a shelling.
 Notice that $<I>_f=0$ is $(I,f)$ if  {shellable}.

\begin{lemma}\label{shel1}

Let  $f:I\mapsto [1,m]$ be such that $(I,f)$ is $\ell$--{ shellable}. Then for every $m_0\le |I|$, there is a subset $S\subset I$
such that $(S,f)$ is $\ell$--{shellable} (and hence $<S>_f=0$ )and $m_0-\ell +1\le |S|\le m_0$.
\end{lemma}
\begin{proof}
Take  a shelling  $I=I_1\cup \cdots \cup I_t$ for $(I,f)$.
Take a maximal $s\le t$ such that $|I_1\cup \cdots \cup I_s|\le m_0$ and put $S=I_1\cup \cdots \cup I_s$.
We must have $|S|\ge m_0-\ell +1$, since otherwise $s\le t-1$, and hence $|I_1\cup \cdots \cup I_{s+1}|\le m_0$, contradicting the
maximality of $s$.
\end{proof}

\begin{proofof}{Theorem}{\ref{weightedgao}}

\noindent{\bf Claim }1. For every any two subsets $A,B\subset [1,m]$, such that $|B|\ge |A|$ there is $R\subset A$
 and an injection $g:R\mapsto  B$ such that $(R,g)$ is $D$--shellable and $|R|\ge |A|-D+1$.

Take a maximal subset  $R\subset A$
and an injection $g:R\mapsto B$ such that $(R,g)$ is $D$ shellable (recall that the empty map is injective and
hence $(\emptyset,\emptyset)$ is $\ell$--shellable). We must have $|R|\ge |A|-D+1$. Suppose that
 $|R|\le |A|-D$. Take a $D$--subset of $A'\subset (A\setminus R)$
 and  $D$--subset $B'\subset (B\setminus g(R))$ . Let $h$ be an arbitrary bijection from $R'$ onto $S'$. By the definition
of the Davenport constant there is a nonempty subset $R'\subset A'$ such that $<R'>_h=0$.
We now extend $g$ to $R\cup R'$ by taking $g(x)=h(x)$ for all $x\in R'$. Now $(R\cup R')$
is  $(D,g)$--shellable contradicting the maximality of $R$.

Assume first $D\le \ell$. Since $m\ge n-1$, we may take any subset $S_0$ of $[1,m]$ with cardinality $n-1$.
By Claim 1, there is $S\subset S_0$
with $n-1\ge |S|\ge n-D+1$ and an injection $g:S\mapsto [1,m]$ such that $<S>_g=0$. The result holds in this case.

 Assume now $\ell < D$.
 Thus $m=n-\ell+D-1$.

\noindent {\bf  Claim} 2. There is a subset $T\subset [1,m]$ such that
 $|T|\ge  D-\ell $  and an injection $g:T\mapsto [1,m]$  such that $(T,g)$ is $\ell$--shellable.
 Take a maximal subset  $R\subset [1,m]$
and an injection $g:R\mapsto [1,m]$ such that $(R,g)$ is $\ell$--shellable (recall that the empty map is injective and
hence $(\emptyset,\emptyset)$ is $\ell$--shellable). We must have $|R|\ge D-\ell$, since otherwise
 $m-|R|\ge n-\ell+D-1-(D-\ell -1)\ge n $. By Lemma there a nonempty subset $R'\subset [1,m]\setminus R$ with
 $|R'|\le \ell$ and an injection  $h: R'\mapsto [1,m]\setminus g(R)$ such that $<R'>_h=0$.
We now extend $g$ to $R\cup R'$ by taking $g(x)=h(x)$ for all $x\in R'$. Now $(R\cup R',g)$
is  $\ell$--shellable contradicting the maximality of $R$.

Put $C=[1,m]\setminus T$. By Claim 1, applied to $C$, there are a subset $E\subset C$ and  injective map $f:E\mapsto [1,m]\setminus g(R)$,
such that $n-1\ge |C|\ge |E|\ge |C|-D+1$ and ${<{E}>}_f=0.$
Now $ n-1-|E|\le n-1-|C|+D-1=(m+\ell-D-1)-|C|+1=|T|+\ell-D<|T|.$
By Lemma \ref{word1} there is subset $F\subset T$ such that $<F>_f=0$ and $n-1-|E|\ge |F|\ge n-1-|E|-\ell+1$.
In particular $n-1\ge |S|\ge n-\ell$, where $S=E\cup F$.
We extend $f$ to $S$, by taking $f(i)=g(i)$ for every $i\in F$. Clearly $f$ is injective on $S$ and
$<S>_f=<E>_f+<F>_f=0$.\end{proofof}

\begin{proofof}{Corollary}{\ref{baricenter}}
 Without loss of generality we may assume $x_{m-r+1}=\cdots =x_m$.
By Theorem \ref{weightedgao}, applied to $(x_1-x_m,\cdots ,x_{m-r}-x_m)$, there are a
subset $I\subset [1,m-r]$ and an injection $f:I\mapsto [1,m-r]$ such that
$n-\min (D,\ell) \le |I| \le n-1$ and
  $$\sum_{ i\in I}w_{i}(x_{f(i)}-x_m)=0,$$ and hence
   $$\sum_{ i\in I}w_{i}x_{f(i)}=(\sum_{ i\in I}w_{i})x_m.$$
  Take now an $(n-|I|)$--subset $J\subset [1,m-r]\setminus I$. We may extend $f$ to an injection on $I\cup J$,
  by putting $f(i)=h(i)$ for every $i\in J$, where $h$ is an arbitrary injection $:J\mapsto [m-r+1,m]$.

  We have
  \begin{eqnarray*}
  \sum_{ i\in I\cup J}w_{i}x_{f(i)}&=&\sum_{ i\in I}w_{i}x_{m}+\sum_{ i\in J}w_{i}x_{m}\\
  &=& (\sum_{ i\in I\cup J} w_{i})x_{m}.
  \end{eqnarray*}\end{proofof}

\end{document}